\documentclass{tac}

\usepackage{amsmath}
\usepackage{amssymb}
\usepackage{proof}
\usepackage{fancybox}
\usepackage{enumerate}
\usepackage[all]{xy}

\usepackage{euscript}
\usepackage{stmaryrd}

\newcommand{\oxl}{ \raisebox{0.3ex}{\ensuremath{\,\scriptscriptstyle\langle}}{\ensuremath{\hspace{-0.2ex}\varotimes\,}} }
\newcommand{\oxr}{{\ensuremath{\,\varotimes}{\raisebox{0.3ex}{\ensuremath{\hspace{-0.2ex}\scriptscriptstyle\rangle\,}}}}}
\newcommand{\otl}{ \raisebox{0.3ex}{\ensuremath{\,\scriptscriptstyle\langle}}{\ensuremath{\hspace{-0.2ex}\varoplus\,}} }
\newcommand{\otr}{{\ensuremath{\,\varoplus}{\raisebox{0.3ex}{\ensuremath{\hspace{-0.2ex}\scriptscriptstyle\rangle\,}}}}}

\newcommand{\prooftri}[1]{*\txt{\begin{picture}(0,0)(10,20)
     	\put(0,0){\line(1,2){10}}
	\put(0,0){\line(1,0){20}}
	\put(20,0){\line(-1,2){10}}
	\put(6,5){\ensuremath{#1}}
	\end{picture}}}

\newenvironment{player}{\( \left\{\begin{array}{c}}{\end{array}\right\} \)}
\newenvironment{opponent}{\( \left(\begin{array}{c}}{\end{array}\right) \)}
\newcommand{\size}{\ensuremath{\mathrm{size}}}
\newcommand{\esize}{\ensuremath{\mathrm{esize}}}
\newcommand{\usize}{\ensuremath{\mathrm{usize}}}
\newcommand{\leaves}{\ensuremath{\mathrm{leaves}}}
\newcommand{\profile}{\ensuremath{\mathrm{prof}}}

\newcommand{\cat}[1]{\ensuremath{\mathbf{#1}}}

\newcommand{\ox}{\ensuremath{\varotimes}}
\newcommand{\bigox}{\ensuremath{\bigotimes}}
\newcommand{\ot}{\ensuremath{\varoplus}}
\newcommand{\bigot}{\ensuremath{\bigoplus}}
\newcommand{\vd}{\ensuremath{\vdash}}

\newcommand{\meet}{\ensuremath{\wedge}}
\newcommand{\join}{\ensuremath{\vee}}
\newcommand{\bigmeet}{\ensuremath{\bigwedge}}
\newcommand{\bigjoin}{\ensuremath{\bigvee}}

\newcommand{\iso}{\ensuremath{\cong}}

\newcommand{\Lra}{\ensuremath{\Longrightarrow}}

\newcommand{\ora}[1]{\ensuremath{\overrightarrow {#1}}}
\newcommand{\ola}[1]{\ensuremath{\overleftarrow {#1}}}
\newcommand{\ol}[1]{\ensuremath{\overline {#1}}}

\newcommand{\wh}[1]{\ensuremath{\widehat{#1}}}

\newcommand{\qqquad}{\quad \quad \quad}

\newcommand{\senta}[1]{\begin{center} #1 \end{center}}

\makeatletter

\newdimen\w@dth

\def\setw@dth#1#2{\setbox\z@\hbox{\scriptsize $#1$}\w@dth=\wd\z@
\setbox\@ne\hbox{\scriptsize $#2$}\ifnum\w@dth<\wd\@ne \w@dth=\wd\@ne \fi
\advance\w@dth by 1.2em}

\def\t@^#1_#2{\allowbreak\def\n@one{#1}\def\n@two{#2}\mathrel
{\setw@dth{#1}{#2}
\mathop{\hbox to \w@dth{\rightarrowfill}}\limits
\ifx\n@one\empty\else ^{\box\z@}\fi
\ifx\n@two\empty\else _{\box\@ne}\fi}}
\def\t@@^#1{\@ifnextchar_ {\t@^{#1}}{\t@^{#1}_{}}}

\def\t@left^#1_#2{\def\n@one{#1}\def\n@two{#2}\mathrel{\setw@dth{#1}{#2}
\mathop{\hbox to \w@dth{\leftarrowfill}}\limits
\ifx\n@one\empty\else ^{\box\z@}\fi
\ifx\n@two\empty\else _{\box\@ne}\fi}}
\def\t@@left^#1{\@ifnextchar_ {\t@left^{#1}}{\t@left^{#1}_{}}}

\def\two@^#1_#2{\def\n@one{#1}\def\n@two{#2}\mathrel{\setw@dth{#1}{#2}
\mathop{\vcenter{\hbox to \w@dth{\rightarrowfill}\kern-4ex
                 \hbox to \w@dth{\rightarrowfill}}%
       }\limits
\ifx\n@one\empty\else ^{\box\z@}\fi
\ifx\n@two\empty\else _{\box\@ne}\fi}}
\def\tw@@^#1{\@ifnextchar_ {\two@^{#1}}{\two@^{#1}_{}}}

\def\tofr@^#1_#2{\def\n@one{#1}\def\n@two{#2}\mathrel{\setw@dth{#1}{#2}
\mathop{\vcenter{\hbox to \w@dth{\rightarrowfill}\kern-4ex
                 \hbox to \w@dth{\leftarrowfill}}%
       }\limits
\ifx\n@one\empty\else ^{\box\z@}\fi
\ifx\n@two\empty\else _{\box\@ne}\fi}}
\def\t@fr@^#1{\@ifnextchar_ {\tofr@^{#1}}{\tofr@^{#1}_{}}}

\newdimen\W@dth
\def\setW@dth#1#2{\setbox\z@\hbox{$#1$}\W@dth=\wd\z@
\setbox\@ne\hbox{$#2$}\ifnum\W@dth<\wd\@ne \W@dth=\wd\@ne \fi
\advance\W@dth by 1.2em}

\def\T@^#1_#2{\allowbreak\def\N@one{#1}\def\N@two{#2}\mathrel
{\setW@dth{#1}{#2}
\mathop{\hbox to \W@dth{\rightarrowfill}}\limits
\ifx\N@one\empty\else ^{\box\z@}\fi
\ifx\N@two\empty\else _{\box\@ne}\fi}}
\def\T@@^#1{\@ifnextchar_ {\T@^{#1}}{\T@^{#1}_{}}}

\def\T@left^#1_#2{\def\N@one{#1}\def\N@two{#2}\mathrel{\setW@dth{#1}{#2}
\mathop{\hbox to \W@dth{\leftarrowfill}}\limits
\ifx\N@one\empty\else ^{\box\z@}\fi
\ifx\N@two\empty\else _{\box\@ne}\fi}}
\def\T@@left^#1{\@ifnextchar_ {\T@left^{#1}}{\T@left^{#1}_{}}}

\def\Tofr@^#1_#2{\def\N@one{#1}\def\N@two{#2}\mathrel{\setW@dth{#1}{#2}
\mathop{\vcenter{\hbox to \W@dth{\rightarrowfill}\kern-4ex
                 \hbox to \W@dth{\leftarrowfill}}%
       }\limits
\ifx\N@one\empty\else ^{\box\z@}\fi
\ifx\N@two\empty\else _{\box\@ne}\fi}}
\def\T@fr@^#1{\@ifnextchar_ {\Tofr@^{#1}}{\Tofr@^{#1}_{}}}

\def\Two@^#1_#2{\def\N@one{#1}\def\N@two{#2}\mathrel{\setW@dth{#1}{#2}
\mathop{\vcenter{\hbox to \W@dth{\rightarrowfill}\kern-4ex
                 \hbox to \W@dth{\rightarrowfill}}%
       }\limits
\ifx\N@one\empty\else ^{\box\z@}\fi
\ifx\N@two\empty\else _{\box\@ne}\fi}}
\def\Tw@@^#1{\@ifnextchar_ {\Two@^{#1}}{\Two@^{#1}_{}}}

\def\to{\@ifnextchar^ {\t@@}{\t@@^{}}}
\def\from{\@ifnextchar^ {\t@@left}{\t@@left^{}}}
\def\two{\@ifnextchar^ {\tw@@}{\tw@@^{}}}
\def\tofro{\@ifnextchar^ {\t@fr@}{\t@fr@^{}}}
\def\To{\@ifnextchar^ {\T@@}{\T@@^{}}}
\def\From{\@ifnextchar^ {\T@@left}{\T@@left^{}}}
\def\Two{\@ifnextchar^ {\Tw@@}{\Tw@@^{}}}
\def\Tofro{\@ifnextchar^ {\T@fr@}{\T@fr@^{}}}

\makeatother

\newcommand{\type}{\ensuremath{::}}
\newcommand{\semi}{\ensuremath{\,;}}
\newcommand{\cn}{\ensuremath{\! : \!}}
\newcommand{\cost}{\ensuremath{\mathrm{cost}}}
\renewcommand{\choose}[2]{\ensuremath{\genfrac{(}{)}{0ex}{1}{#1}{#2}}}

\title{On the Complexity of Cockett-Seely Polarized Games}
\author{J. R. B. Cockett and C. A. Pastro}
\thanks{Research partially supported by NSERC, Canada. Diagrams were produced
with the \Xy-pic package of K. Rose and R. Moore and inferences with
M. Tatsuya's \texttt{proof.sty}.}
\address{Department of Computer Science, University of Calgary \\
2500 University Drive NW, Calgary, AB, Canada T2N 1N4}
\eaddress{robin@cpsc.ucalgary.ca \quad  pastroc@cpsc.ucalgary.ca}
\copyrightyear{2004}

\begin{document}
\maketitle

\begin{abstract}
In this paper the complexity of provability of polarized additive,
multiplicative, and exponential formulas in the (initial) Cockett-Seely
polarized game logic is discussed. The complexity is ultimately based on the
complexity of finding a strategy in a formula which is, for polarized
additive formulas, in the worst case linear in their size. Having a proof
of a sequent is equivalent to having a strategy for the 
internal-hom object. In order to show that the internal-hom object can have
size exponentially larger than the formulas of the original sequent we
develop techniques for calculating the size of the multiplicative formulas.

The structure of the internal hom object can be exploited and, using dynamic 
programming techniques, one can reduce the cost of finding a strategy in
such a formula to the order of the product of the sizes of the original
formulas. The use of dynamic techniques motivates the consideration of games
as acyclic graphs and we show how to calculate the size of these graph games
for the multiplicative and additive fragment and, thus, the cost of
determining their provability using this dynamic programming approach.

The final section of the paper points out that, despite the apparent
complexity of the formulas, there is, for the initial polarized logic with
all the connectives (additives, multiplicatives, and exponentials) a way of
determining provability which is \emph{linear} in the size of the formulas.
\end{abstract}

\section*{Introduction}

It is a natural question in a logic to ask how hard is it to decide whether a
sequent is provable. This paper examines this question for the Cockett-Seely
polarized game logic \cite{cockett04:polarized} which we will henceforth refer 
to simply as polarized game logic. We shall show that this question can be
answered in linear time.

The decision procedure consists of translating the given sequent in the
polarized game logic into a single polarized game using ``negation'' and the
``mixed par'' operation. As this is an internal-hom game, determining whether
a strategy exists in this game corresponds to determining whether a proof
exists for the sequent. 

To determine whether a strategy exists in a game requires, in the worst
case, time proportional to the number of edges in the tree representing
that game.  Thus, it becomes important to estimate the size of this
internal-hom game.  While, in general, there is no simple expression for
this size, it is often possible to estimate the ``uniform size'' of an
internal-hom game. The uniform size is defined as the number of binary
logical operations required, in the worst case, to determine whether a
strategy is present. This makes it possible to get a fairly precise
estimate for the size of certain internal-hom games. These estimates
allow us to demonstrate that there can be an exponential increase in the
size of the internal-hom game over the original games.  

Fortunately it is possible to exploit the structure of the internal hom 
object and, using dynamic programming techniques, reduce the cost of finding
a strategy in such an object to be bounded by the product of the sizes of
the original game trees. This approach involves regarding games as if they
are acyclic directed graphs, a direction already pioneered by Hyland and
Schalk~\cite{hyland-schalk}. When this is done one can give precise
expressions for the size of formulas built using the multiplicatives
connectives and, thus, for the cost of the dynamic programming approach to
finding a strategy.

The categorical semantics of polarized games is outside the scope of this
paper. In order to keep the paper largely self-contained we do, however,
provide a description and brief discussion of the basic concepts. The reader
who is interested in a more complete story should consult the paper of
Cockett and Seely~\cite{cockett04:polarized}. Such a course of action is
recommended as, in the last section of this paper, we use this semantics,
albeit in a very simple way, to show that for this model one can decide
provability in linear time for the full logic with polarized additives,
multiplicatives, and exponentials.

Before starting this story we should emphasize that the multiplicative
structure of this model is \emph{not} the free structure. This model has
``soft'' multiplicative and exponential structures
(see~\cite{cockett04:polarized}). Provability in free polarized
multiplicative structures has the same complexity as provability in the 
corresponding multiplicative fragment of linear logic (as any multiplicative
proof net can be polarized).  For multiplicative linear logic with atoms
provability is known to be NP-complete \cite{kanovich92:complexity}, and
even without atoms this decision problem is just as hard: the multiplicative
units can be used to simulate the effect of having
atoms~\cite{lw94:complexity}. In the presence of atoms, additives, and
exponentials (or even nonlogical axioms) these provability problems become
undecidable for ordinary linear logic~\cite{lmss92:complexity}.  

The results we obtain are somewhat curious as they apply to the finite
fragment of precisely the same model whose depolarization was used by
Abramsky and Jagadeesan~\cite{abramsky94:games} to obtain full completeness
for the (iso-mix) multiplicative fragment of multiplicative linear logic.
This seems to highlight the important role units play in providing these
logics with their underlying complexity. Note also that using the
transformation to Laurent's version of polarized
logic~\cite{laurent02:polarized}, discussed in~\cite{cockett04:polarized},
this model also provides, through the coKleisli construction for the $!(\_)$
comonad, a non-trivial model of intuitionistic logic which has a linear
time provability. It is worth recalling that this is, for the free logic,
a P-space complete problem~\cite{statman}.

The outline of this paper is as follows. In Section~\ref{sec-polgam}
polarized games are introduced. In Section~\ref{sec-pollog} polarized game
logic is introduced. Section~\ref{sec-mult} describes the multiplicative
and exponential structure on polarized games. Section~\ref{sec-strat} shows
how provability can be turned into the question of finding a strategy in the
internal-hom object. Section~\ref{sec-size} discusses the size of the
multiplicative formulas including internal-hom objects. In
Section~\ref{sec-dyn}, we discuss the relationship between the dynamic
programming solution to determining provability and graph games. In
particular, we describe simply formulas for the size of multiplicatives in
graph games. Finally, in Section~\ref{sec-prov}, we show that provability,
despite the evident complexity of the underlying additive structures, can
be evaluated for this model in linear time for the complete additive,
multiplicative, and exponential logic.

\section{Polarized games} \label{sec-polgam}

A polarized game is a finite 2-player input-output game of the sort studied
by Abramsky and Jagadeesan~\cite{abramsky94:games} and Hyland
in~\cite{hyland97:games}. The two players are referred to as the ``opponent''
and the ``player''. Games are often discussed using the terminology of
processes; in these terms a game is a protocol for interaction while a
map between games is a process which interacts through these protocols on
the ``input'' (or domain) and ``output'' (or codomain) channels. On the
output channel, the opponent messages come from the ``environment'' and the
player messages are generated by the process or ``system''. On the input
channel, however, the roles are reversed: the opponent messages are
generated by the ``system'' and the player messages come from the
``environment.''  

There are two sorts of games: those in which the opponent has the first move
and those in which the player has the first move. The games introduced in
this paper are the games of the initial polarized game
category~\cite{cockett04:polarized}. Games are represented in a number of
ways.

An \textbf{opponent game} is denoted by
\[
O = (b_1 \cn P_1,\ldots,b_m \cn P_m) = \bigsqcap\limits_{i \in I} b_i \cn P_i =
\vcenter{\xymatrix@R=2ex@C=3ex@M=0ex{
& \circ \ar@{-}[dl]_{b_1} \ar@{-}[d] \ar@{-}[dr]^{b_m} & \\ P_1 & \cdots & P_m}}
\]
where $I = \{1,\ldots,m\}$ and each $P_i$ is, inductively, a player game.
A \textbf{player game} is denoted by
\[
P = \{a_1 \cn O_1,\ldots,a_n \cn O_n\} = \bigsqcup\limits_{j \in J} a_j \cn O_j=
\vcenter{\xymatrix@R=2ex@C=3ex@M=0ex{
& \bullet \ar@{-}[dl]_{a_1} \ar@{-}[d] \ar@{-}[dr]^{a_n} & \\
O_1 & \cdots & O_n}}
\]
where $J = \{1,\ldots,n\}$ and each $O_j$ is, inductively, an opponent game.
We allow (opponent and player) games with empty index sets; this gives two
atomic games at which the inductive construction of games begins:
\[
\mathbf{1} = (\,) = \bigsqcap\nolimits_\emptyset = \circ
\qqquad 
\mathbf{0} = \{\, \} = \bigsqcup\nolimits_\emptyset = \bullet
\]

\begin{example}
A typical game (without branch labels) looks like:
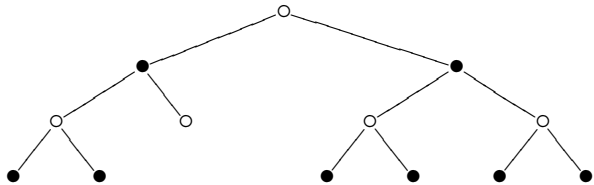
\[\xymatrix@C=2ex@R=3ex@M=0ex{ 
&&&&&&& \circ \ar@{-}[dllll] \ar@{-}[drrrr] \\
&&& \bullet \ar@{-}[dll] \ar@{-}[dr] &&&&&&&& \bullet \ar@{-}[dll]
\ar@{-}[drr] \\
& \circ \ar@{-}[dl] \ar@{-}[dr] &&& \circ &&&&& \circ \ar@{-}[dl] \ar@{-}[dr]
&&&& \circ \ar@{-}[dl] \ar@{-}[dr] \\
\bullet && \bullet &&&&&& \bullet && \bullet && \bullet && \bullet}
\]
Notice that this is a opponent start game.
\end{example}

The \textbf{dual} of a game $G$, denoted $\ol{G}$, is defined as
\[\ol{O} = \ol{\bigsqcap_{i \in I} b_i \cn P_i} = \bigsqcup_{i \in I} b_i \cn
\ol{P_i}
\qqquad 
\ol{P} = \ol{\bigsqcap_{j \in J} a_j \cn O_j} = \bigsqcup_{j \in J} a_j \cn
\ol{O_j}
\]
Thus, taking the dual of a game simply flips the black and white nodes
while leaving the structure and branch labels (the $b_i$'s and $a_j$'s) alone.

\subsection{Morphisms of games} \label{sec-mor-gam}

There are three types of morphisms: opponent morphisms (from opponent games
to opponent games), mixed (or cross) morphisms (from opponent games to player
games), and player morphisms (from player games to player games).  These are defined 
inductively by:

\begin{description}
\item[Opponent morphisms:] 
\[\left(\begin{array}{c}
b_1 \mapsto h_1 \\
\cdots \\
b_m \mapsto h_m \\
\end{array}\right) 
:\xymatrix{O \ar[r] & (b_1 \cn P_1,\ldots,b_m \cn P_m)}
\]
where each $\xymatrix{h_i : O \ar[r] & P_i}$ is a mixed morphism 
(the process listens on the codomain channel).

\item[Mixed morphisms:] These are either of the form
\[
\ola{b_k} \cdot f : \xymatrix{(b_1 \cn P_1,\ldots,b_m \cn P_m) \ar[r] & P}
\]
where $k \in \{1,\ldots,m\}$ and $\xymatrix{f:P_k \ar[r] & P}$ is a
player morphism (the process outputs $b_k$ on the domain channel) or:
\[
\ora{a_k} \cdot g : \xymatrix{O \ar[r] & \{a_1 \cn O_1,\ldots,a_n \cn O_n\}}
\]
where $k \in \{1,\ldots,n\}$ and $\xymatrix{g:O \ar[r] & O_k}$ is an
opponent morphism (the process outputs $a_k$ on the codomain channel).

\item[Player morphisms:] 
\[\left\{\begin{array}{c}
a_1 \mapsto h_1 \\
\cdots \\
a_n \mapsto h_n \\
\end{array}\right\}
:\xymatrix{\{a_1 \cn O_1,\ldots,a_n \cn O_n\} \ar[r] & P}
\]
where each $\xymatrix{h_i \cn O_i \ar[r] & P}$ is a mixed morphism 
(the process listens on the domain channel).

\end{description}
There are no morphisms from player to opponent.

Viewed as a process, a mixed morphisms initiates the interaction by sending
a message down either the domain channel or the codomain channel. Having sent
a message down one of these channels the process must then wait until it
receives a reply from the environment along the same channel. Having received
this reply the process is, once again, in a mixed state and can send a message
down either channel. This process discipline which requires that the reply on
a channel must be received before any further action is permitted means that
processes are completely sequential. In fact, an interacting system of such
processes is ``polarized'' in the very strong sense that at any given time
only one processes can ever be active (i.e., between receiving and sending a
message).

In Hyland's ``games for fun'' \cite{hyland97:games} a process is permitted to
simply not respond when it is in a mixed state, in other words partial
processes are allowed. The processes we consider here must produce a message
whenever it is in a mixed state so that, for us, processes are total. These
correspond, for finite games, to winning strategies in the sense 
of~\cite{abramsky94:games}.

\begin{example} \label{exam-oppmap}
The following is a map between two opponent games:
\[\left(\begin{array}{l}
a \mapsto \ora{c} \cdot (\,) \\
b \mapsto \ora{e} \cdot (\,)
\end{array} \right):\
\vcenter{\xymatrix@C=2ex@R=3ex@M=0ex{
&\circ \ar@{-}[dl]_a \ar@{-}[dr]^b \\
\bullet && \bullet}}
\xymatrix{\ar[r] &}
\vcenter{\xymatrix@C=2ex@R=3ex@M=0ex{
&&& \circ \ar@{-}[dll]_a \ar@{-}[drr]^b \\
& \bullet \ar@{-}[dl]_c \ar@{-}[dr]^d &&&& \bullet\ar@{-}[dl]_e \ar@{-}[dr]^f \\
\circ && \circ && \circ && \circ}}
\]

Notice that this is not the only possible map between these two games, 
the following are also maps:
\[\left(\begin{array}{l}
a \mapsto \ora{d} \cdot (\,) \\
b \mapsto \ora{f} \cdot (\,)
\end{array} \right),
\qquad
\left(\begin{array}{l}
a \mapsto \ola{a} \cdot \{\, \} \\
b \mapsto \ola{b} \cdot \{\, \}
\end{array} \right),
\qquad
\left(\begin{array}{l}
a \mapsto \ora{c} \cdot (\,) \\
b \mapsto \ola{a} \cdot \{\, \}
\end{array} \right)
\]
\end{example}

\subsection{Composition} \label{sec-comp}

There are four types of composition we must consider: opponent opponent
composition, opponent mixed composition, mixed player composition, and player
player composition. These compositions can be defined via rewriting rules 
as follows:

\begin{description}
\item[Opponent opponent composition:]
\[
f \semi 
\left(\begin{array}{c}
b_1 \mapsto h_1 \\
\cdots \\
b_n \mapsto h_n
\end{array}\right)
\Lra
\left(\begin{array}{c}
b_1 \mapsto f \semi h_1 \\
\cdots \\
b_n \mapsto f \semi h_n
\end{array}\right)
\]

\item[Opponent mixed composition:]
\[
\left(\begin{array}{c}
b_1 \mapsto h_1 \\
\cdots \\
b_n \mapsto h_n
\end{array}\right)
; \ola{b_k} \cdot v \ \Lra\ h_k \semi v
\quad \text{and} \quad
f \semi \ora{a} \cdot v \ \Lra\ \ora{a} \cdot (f \semi v)
\]

\item[Mixed player composition:]
\[
\ora{a_k} \cdot v \semi
\left\{\begin{array}{c}
a_1 \mapsto h_1 \\
\cdots \\
b_n \mapsto h_n
\end{array}\right\} \ \Lra\ v \semi h_i
\quad \text{and} \quad
\ola{b} \cdot v \semi g \ \Lra\ \ola{b} \cdot (v \semi g)
\]

\item[Player player composition:]
\[
\left\{\begin{array}{c}
a_1 \mapsto h_1 \\
\cdots \\
a_n \mapsto h_n
\end{array}\right\} ; g \Lra
\left\{\begin{array}{c}
a_1 \mapsto h_1 \semi g \\
\cdots \\
a_n \mapsto h_n \semi g
\end{array}\right\}
\]
\end{description} 

\begin{example}
Here we compose an opponent map (from Example~\ref{exam-oppmap}) with a mixed
map.
\[
\vcenter{\xymatrix@C=2ex@R=3ex@M=0ex{
&\circ \ar@{-}[dl]_a \ar@{-}[dr]^b \\
\bullet && \bullet}}
\xymatrix{\ar[r] &}
\vcenter{\xymatrix@C=2ex@R=3ex@M=0ex{
&&& \circ \ar@{-}[dll]_a \ar@{-}[drr]^b \\
& \bullet \ar@{-}[dl]_c \ar@{-}[dr]^d &&&& \bullet\ar@{-}[dl]_e \ar@{-}[dr]^f \\
\circ && \circ && \circ && \circ}}
\xymatrix{\ar[r] &}
\vcenter{\xymatrix@C=2ex@R=3ex@M=0ex{
&\bullet \ar@{-}[dl]_a \ar@{-}[dr]^b \\
\circ && \circ}}
\]
\[
\begin{array}{rcl}
\left(\begin{array}{l}
a \mapsto \ora{c} \cdot (\,) \\
b \mapsto \ora{e} \cdot (\,)
\end{array} \right)\ ;\ 
\ola{b} \cdot
\left\{\begin{array}{l}
e \mapsto \ora{a} \cdot (\,) \\
f \mapsto \ora{b} \cdot (\,)
\end{array} \right\}
& \Lra &
\ora{e} \cdot (\,) \ ;\ 
\left\{\begin{array}{l}
e \mapsto \ora{a} \cdot (\,) \\
f \mapsto \ora{b} \cdot (\,)
\end{array} \right\} \smallskip\\
& \Lra &
(\,) \ ;\ 
\ora{a} \cdot (\,) \smallskip\\
& \Lra &
 
\ora{a} \cdot ((\,) \semi (\,)) \smallskip\\
& \Lra &
\ora{a} \cdot (\,) 
\end{array}
\]
\end{example}

Composition, as it is defined here, is exactly the cut elimination rewriting
for the polarized game logic, which we introduce in the next section. It can
be shown by an inductive argument that this is a confluent and terminating
rewriting which eliminates the composition. (See the paper of Cockett and
Seely \cite{cockett01:finite} where this is proved for a similar system.)

\begin{proposition} \quad
\begin{enumerate}[{\upshape (i)}]
\item The above rewriting on morphisms terminates.
\item The above rewriting is confluent.
\item The associative law is satisfied by all composite triples.
\end{enumerate}
\end{proposition}

Given a opponent game $O = (b_1 \cn P_1,\ldots,b_m \cn P_m)$, the
\textbf{identity opponent morphism} $\xymatrix{1_O:O \ar[r] & O}$ is defined
inductively as:
\[1_O =
\left(\begin{array}{c}
b_1 \mapsto \ola{b_1} \cdot 1_{P_1} \\
\cdots \\
b_m \mapsto \ola{b_m} \cdot 1_{P_m}
\end{array}\right)
\]

Given a player game $P = \{a_1 \cn O_1,\ldots,a_n \cn O_n\}$, the
\textbf{identity player morphism} $\xymatrix{1_P:P \ar[r] & P}$ is defined
inductively as:
\[1_P =
\left\{\begin{array}{c}
a_1 \mapsto \ora{a_1} \cdot 1_{O_1} \\
\cdots \\
a_n \mapsto \ora{a_n} \cdot 1_{O_n}
\end{array}\right\}
\]

Observe that in any possible composition, the identity acts as a neutral
element, that is $f \semi 1 = f = 1 \semi f$.  This gives:

\begin{proposition} \quad
\begin{enumerate}[{\upshape (i)}]
\item Player games and player morphisms form a category.
\item Opponent games and opponent morphisms form a category.
\item Mixed morphisms between the opponent category and the player category
form a module between these categories.
\end{enumerate}
\end{proposition}

The \textbf{polarized categories} of~\cite{cockett04:polarized} are precisely 
modules $\wh{\cat{C}}:\cat{C}_o \to \cat{C}_p$ and in that paper the details
of how these provide the categorical semantics for the proof theory of
polarized game logic which we now present is spelled out.

\section{Polarized game logic} \label{sec-pollog}

Proofs in polarized game logic with no atomic types or non-logical axioms
correspond precisely to morphisms between polarized games as described above.
Below we display this logic as a Gentzen style sequent calculus. It requires
three kinds of sequent:

\begin{description}
\item[Opponent sequents:] $O \vd_o O'$ where $O$ and $O'$ are opponent 
propositions.

\item[Mixed sequents:] $O \vd P$ where $O$ is an opponent proposition 
and $P$ is a player proposition. (This notation differs slightly from
the notation in~\cite{cockett04:polarized} in that we leave the mixed sequent
turnstile unannotated.)

\item[Player sequents:] $P \vd_p P'$ where $P$ and $P'$ are player 
propositions.
\end{description}

The valid inferences of the polarized game logic are as follows, which are a
``graded'' version of $\Sigma\Pi$~\cite{cockett01:finite}.

\begin{center}
\ovalbox{
\parbox[c]{52ex}{\begin{center}
\medskip
$\infer{A \vd_p A}{} \quad \text{(atomic identities)} \quad
\infer{B \vd_o B}{}$ 
\bigskip\\
$\infer[\text{(cotuple)}]{\bigsqcup\limits_{i \in I} X_i \vd_p Y}{\{X_i \vd Y 
\}_{i \in I}} \qquad
\infer[\text{(tuple)}]{X \vd_o \bigsqcap\limits_{i \in I} Y_i}{\{X \vd Y_i
\}_{i \in I}}$
\bigskip\\
$\infer[\text{(projection)}]{\bigsqcap\limits_{i \in I} X_i \vd Y}{X_k \vd_p Y}
\quad
\infer[\text{(injection)}]{X \vd \bigsqcup\limits_{i \in I} Y_i}{X \vd_o Y_k}$
\medskip
\end{center}}}
\end{center}

Notice that in the cotuple and tuple rules the index set $I$ may be empty, 
though not in the injection and projection rules. Also, observe that the
inference system is symmetric, that is, it has an obvious
$\bigsqcup-\bigsqcap$ symmetry. Explicitly, we may flip the direction of
the sequents while swapping $\bigsqcup$ for $\bigsqcap$ and ``opponent'' for
``player'' to obtain the same system. This symmetry arises from an underlying
categorical duality.

The logic has four cut rules which correspond to those permitted by the
types. The first two arise as cuts in the player and opponent sequents. The
last two are the two possible cuts on the mixed sequent. These correspond
categorically to the actions expected of a module.

\begin{center}
\ovalbox{
\parbox[c]{55ex}{\begin{center}
\medskip
$\infer[\text{(p-cut)}]{X \vd_p Z}{X \vd_p Y & Y \vd_p Z} \quad\
\infer[\text{(o-cut)}]{X \vd_o Z}{X \vd_o Y & Y \vd_o Z}$
\bigskip\\
$\infer[\text{(mp-cut)}]{X \vd Z}{X \vd Y & Y \vd_p Z} \quad
\infer[\text{(om-cut)}]{X \vd Z}{X \vd_o Y & Y \vd Z}$
\medskip
\end{center}}}
\end{center}

It is an easy exercise to prove:

\begin{theorem}
The polarized game logic satisfies cut elimination. 
\end{theorem}

In fact, cut elimination is exactly the rewriting steps introduced in the
previous section for composition; we can make this explicit by annotating 
the logic with terms. Following \cite{cockett04:polarized} we use ``$\type$''
to denote the term-type membership relation, e.g., $t \type U \vd V$ will
mean that $t$ is a term of type $U \vd V$, where $U$ (say) may be of the form
$a \cn X$.

\begin{center}
\ovalbox{
\parbox[c]{80ex}{\begin{center}
\medskip
$\infer{1_A \type A \vd_p A}{} \quad \text{(atomic identities)} \quad
\infer{1_B \type B \vd_o B}{}$
\bigskip\\
$\infer[\text{(cotuple)}]{\{a_i \cn h_i\}_{i \in I} \type \bigsqcup\limits_{i
\in I} a_i \cn X_i \vd_p Y}{\{h_i \type X_i \vd Y\}_{i \in I}} \quad
\infer[\text{(tuple)}]{(b_i \cn h_i)_{i \in I} \type X \vd_o 
\bigsqcap\limits_{i \in I} b_i \cn Y_i}{\{h_i \type X \vd Y_i \}_{i \in I}}$
\bigskip\\
$\infer[\text{(projection)}]{\ola{a_i} \cdot g \type \bigsqcap\limits_{i \in I}
a_i \cn X_i \vd Y}{g \type X_k \vd_p Y} \quad
\infer[\text{(injection)}]{\ora{b_k}\cdot f \type X \vd \bigsqcup\limits_{i
\in I} b_i \cn Y_i}{f \type X \vd_o Y_k}$
\bigskip\\
$\infer[\text{(cut)}]{f \semi g \type X \models Z}{f \type X \models Y & g
\type Y \models Z}$ 
\smallskip\\
where $\models$ represents the appropriate type \\
of sequent for each of the four cut rules
\medskip
\end{center}}}
\end{center}

As this logic satisfies cut elimination, to determine provability of a
sequent it suffices to establish that there is a cut free proof of the
sequent. This is precisely the same as establishing that there is a morphism
of games which is in normal form (i.e., one in which the composition (cut)
has been removed).

\section{Multiplicative and exponential structure} \label{sec-mult}

In this section we describe the multiplicative structure on polarized games.
This structure will be used in the calculation of the complexity of provability. 
We shall also introduce the exponential types (the Curien exponential) although 
we shall not develop their properties in any depth.

Suppose we have games:
\[O = (b_i \cn P_i \mid i \in I),\ O' = (b_i' \cn P_i' \mid i \in I'),\
P = \{a_j \cn O_j \mid j \in J\},\ P' = \{a_j' \cn O_j' \mid j \in J'\}
\]

\begin{description}
\item[{$[O \oxr P]$}] This operation takes an opponent game $O$ and a player
game $P$ and produces a player game. The operation is defined inductively
by:
\[O \oxr P = \{a_j \cn O \ox O_j \mid j \in J\}
\]
where the tensor operation $\ox$ is defined below. Notice that the
``direction'' of the operation points to the game whose type is inherited
by $O \oxr P$.

\item[{$[P \oxl O]$}] This operation is a simple left-right dual of the
preceding one, taking a player game $P$ and an opponent game $O$ and 
producing a player game. The operation is defined inductively by:
\[P \oxl O = \{a_j \cn O_j \ox O \mid j \in J\}
\]

\item[{$[O \ox O']$}] This operation takes in two opponent games and produces
an opponent game. It is defined inductively by:
\[O \ox O' = (b_i \cn P_i \oxl O', b_k' \cn O \oxr P_k' \mid i \in I, k \in I')
\]

\item[{$[P \otr O]$}] This operation takes a player game $P$ and an opponent
game $O$ and produces an opponent game. The operation is defined inductively
by:
\[P \otr O = (b_i \cn P \ot P_i \mid i \in I)
\]
where the par operation $\ot$ is defined below. Again, the direction of
operation points to the game whose type is inherited by the compound game.

\item[{$[O \otl P]$}] This operation is the left-right dual, taking an
opponent game $O$ and a player game $P$ and producing an opponent game.
The operation is defined inductively by:
\[O \otl P = (b_i \cn P_i \ot P \mid i \in I)
\]

\item[{$[P \ot P']$}] This operation takes in two player games and produces
a player game. It is defined inductively by:
\[P \ot P' = \{a_j \cn O_j \otl P', a_k' \cn P \otr O_k' \mid j \in J, k \in J'\}
\]

\item[{$[!O]$}] The exponential (bang) is defined inductively by:
\[!O = \bigox_{i \in I} (b_i:!'P_i) \text{ for $|I| > 0$}
\quad \text{and} \quad
!(\,) = (\,)
\]
where $!'P = \{a_j: !O_i \mid j \in J\}$.

\item[{$[?P]$}] The exponential (whimper) is defined inductively by:
\[?P = \bigot_{j \in J} \{a_j:?'O_j\} \text{ for $|J| > 0 $}
\quad \text{and} \quad
?\{\,\} = \{\,\}
\]
where $?'O = (b_i: ?P_i \mid i \in I)$.

\end{description}

Notice that there are some obvious dualities:
\[O \ox O = \ol{\ol{O} \ot \ol{O}} \qquad
O \oxr P = \ol{\ol{O} \otl \ol{P}} \qquad
P \oxl O = \ol{\ol{P} \otr \ol{O}} \qquad 
!O = \ol{?\ol{O}}
\] 

There are a number of simple identities that hold of these operations; we
list (without proof) several for reference.

\begin{lemma}
\[\begin{array}{rclcrcl}
O \ox O' & \iso & O' \ox O && P \ot P' & \iso & P' \ot P \\
O \ox \mathbf{1} & \iso &  O && \mathbf{0} \ot P & \iso & P \\
O \oxr \mathbf{0} & \iso & \mathbf{0} && \mathbf{1} \otl P & \iso & \mathbf{1}\\
\mathbf{1} \oxr P & \iso & P && O \otl \mathbf{0} & \iso & O \\
(O \ox O') \oxr P & \iso & O \oxr (O' \oxr P) 
    && O \otl (P \ot P) & \iso & (O \otl P) \otl P') \\
(O_1 \ox O_2) \ox O_3 & \iso & O_1 \ox (O_2 \ox O_3)
    && (P_1 \ot P_2) \ot P_3 & \iso & P_1 \ot (P_2 \ot P_3) \\
\end{array}
\]
\end{lemma}

These combinatoric operations are functors on the game categories.
In the opponent category there is a bijective correspondence
\[\infer={\xymatrix{O_2 \ar[r] & \ol{O_1} \otr O_3}}
{\xymatrix{O_1 \ox O_2 \ar[r] & O_3}}
\]
which makes it into a monoidal closed category, see~\cite{cockett04:polarized}.
The exponential $!(\_)$ gives rise to a comonad in whose Kleisli category 
the tensor turns into a product.

\section{Provability and strategies} \label{sec-strat}

There is a natural question to ask when one regards games as a logic, 
viz, how hard is it to decide whether a sequent is provable. Recall 
that there are three sorts of sequents:
\[O \vd_o O', \quad O \vd P, \quad \text{and} \quad P \vd_p P'.
\]

We start by pointing out that it suffices to determine the complexity of
deciding whether an opponent sequent $O \vd_o O'$ is provable. This is 
because if we are faced with a player sequent we may negate the games and
view it as an opponent sequent, i.e.,
\[P \vd_p P' = \ol{P'} \vd_o \ol{P}
\]
Similarly any problem in the mixed setting $O \vd P$ may be transformed to a
problem in the opponent setting by adding an apr\`es vous move to the player
game. That is, we simply prepend the player game with an opponent node to
turn it into an opponent game, e.g.,
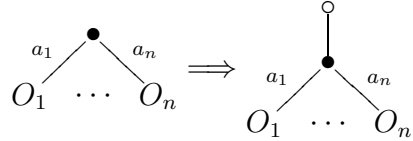
\[\vcenter{\xymatrix@R=3ex@C=2ex@M=0ex{
& \bullet \ar@{-}[dl]_{a_1} \ar@{-}[dr]^{a_n} \\
O_1 & \cdots & O_n}} \Lra
\vcenter{\xymatrix@R=3ex@C=2ex@M=0ex{
& \circ \ar@{-}[d] \\ & \bullet \ar@{-}[dl]_{a_1} \ar@{-}[dr]^{a_n} \\
O_1 & \cdots & O_n}}
\]
Logically this is because:
\[ \infer={O \vd P}{O \vd_o \bigsqcap\nolimits_* P}
\]

Therefore, it suffices to determine how difficult it is to decide if there
is a proof of an opponent sequent $O \vd_o O'$. There is a special case
where this determination is particularly simple, namely when the game on the
left-hand side of the turnstile is $\circ = (\,)$. In this case we are simply
asking whether there is a morphism $\xymatrix{f:(\,) \ar[r] & O'}$, which we
shall call a \textbf{strategy} in the opponent game $O'$.

The calculation below will determine whether or not a strategy exists
in a game $G$.  Notice that the same procedure can be used to count the
number of strategies by turning \textsf{true} to 1, \textsf{false} to 0,
meets to products, and joins to sums. The calculation is as follows:

\begin{itemize}
\item The empty opponent and the empty player games 
respectively may be evaluated respectively as having a strategy and 
not having a strategy respectively, i.e.,
\senta{$\left\| (\,) \right\| = \mathsf{true} \qquad
\left\| \{\,\} \right\| = \mathsf{false}$}

\item At an opponent node $\bigsqcap_{i \in I} P_i$ we must find a strategy
in each of the $P_i$, and so there is a strategy for $\bigsqcap_{i \in I} P_i$
if there is a strategy in the meet of the $P_i$'s. This is indicated by:
\[\left\| \bigsqcap\nolimits_{i \in I} P_i \right\| =
\bigmeet\nolimits_{i \in I} \| P_i \|
\]

\item At a player node $\bigsqcup_{i \in I} O_i$ we must find a strategy in
any one of the $O_i$'s, and so there is a strategy for $\bigsqcup_{i \in I}
O_i$ if there is a there is a strategy in the join of the $O_i$'s. This is
indicated by:
\[\left\| \bigsqcup\nolimits_{i \in I} O_i \right\| =
\bigjoin\nolimits_{i \in I} \| O_i \|
\]
\end{itemize}

Therefore we obtain a decision procedure which is, in the worst case, linear
in the size of the game.  

\begin{remark} \label{rem-strat}
\begin{enumerate}
\item A \textbf{counter-strategy} corresponds to a proof of a sequent
$P \vd_p \{\,\}$ or $O \vd \{\,\}$. Notice that the determination of whether
there is a counter-strategy is exactly the (de Morgan) dual calculation in
which \textsf{true} and \textsf{false} are swapped as are $\join$ and
$\meet$. This has the consequence that a strategy exists for a game $G$ if
and only if $G$ has no counter-strategy.

\item The caveat that this is a worst case scenario is important as in most
cases one does not have to visit the whole tree. Indeed, the more ``bushy''
the tree the more likely this is!  We may do the following very rough
calculation to determine the expected time it takes to determine whether a
(bushy) tree of maximum depth $d$ has a strategy. We suppose the probability
that an (opponent) game tree has a strategy is $p$. This means the
probability that the dual (player) tree has a counter-strategy is also $p$.
At an opponent node we are looking for the first player subtree which does
not have a strategy (i.e., has a counter-strategy) as then we can cut off
the evaluation of the conjunct and return false. Dually at a player node we
are looking for the first opponent subtree which has a strategy as then we
can cut off the evaluation of the disjunct and return true.  

If we are sitting at an opponent node with a very large branching factor the
expected number of edges we will have to inspect is then the sum of the
probabilities that the edge will have to be inspected. This is therefore
bounded by:
\[\sum_{i=0}^\infty p^i = \frac{1}{1-p}
\]
Now this happens at each level of the tree so that the expected cost of
evaluating a tree is $(\frac{1}{1-p})^d$, where $d$ is the depth of the tree.  

Thus, rather surprisingly, the fact that a tree is bushy does not affect the 
expected evaluation time as adversely as might have been thought. A
reasonable guess at the probability $p$ is, of course, $\frac{1}{2}$ so that
to evaluate a tree of depth $d$ might be expected to take time $2^d$ no
matter how bushy it is!
\end{enumerate}
\end{remark}

We now need to be more precise about our calculations and, in particular,
what we mean by ``size'' of a game tree $G$. We will use a variety of
measures. The first and simplest is the number of edges in the game tree,
which we write $\esize[G]$. Note that the edges may be counted first by
those which leave opponent nodes, $\esize_o[G]$, and then those which leave
player nodes, $\esize_p[G]$. Very closely related to this measure (for
trees) is the number of nodes $|G|$, which we shall mean when we speak of
the size of $G$. Once again the nodes can be divided into the number of
opponent nodes, $|G|_o$, and the number of player nodes, $|G|_p$. It is not
hard to see that $\esize[G]=|G|-1$ for any game tree.

Another way of measuring the size of a game tree is the \textbf{uniform size}
which is defined as:
\begin{itemize}
\item $\usize\left[(\,)\right] = \usize\left[\{\,\}\right] = 0$
\item $\usize\left[\bigsqcap_{i \in I} a_i \cn G_i\right] = \usize\left[
\bigsqcup_{i \in I} a_i \cn G_i\right] = n-1 + \sum_{i \in I} \usize[G_i]$ 
\end{itemize}
The form of this definition immediately tells us that $\usize[G] \leq
\esize[G]$. The uniform size of a game tree is important as it measures the
number of binary logical operations required to evaluate the tree. To see
this note that to evaluate an $m$-ary conjunction requires $m-1$ binary
operations.

Yet another way to measure the size of a game tree is by counting its leaves,
which is written $\leaves[G]$. Another reason why the uniform size is a
simple measure is that it is exactly the number of leaves less one. We may
summarize the relationship between all these size measurements as follows:
\[\usize[G] = \leaves[G]-1 \leq \esize[G] = |G|-1 = \esize_o[G]+\esize_p[G]
= |G|_o+|G|_p-1
\]

The difference between $\usize[G]$ and $\esize[G]$ may be viewed as the cost
associated with retrieving the arguments of the binary operations.
Notice that the uniform cost of any tree with branching factor one is
therefore zero as the complete cost of evaluation is in the retrieval! In
fact, it is precisely nodes with branching factor one which cannot be
accounted for (upto a constant factor) by the uniform size. 
In fact, if one knows the depth of the game tree is say $d$, then each path
to a leaf contains at most $d$ edges and so one can bound the edge size by:
\[\usize[G] \leq \esize[g] \leq d \cdot (1+ \usize[G])
\]
This means that the uniform size of a game tree is generally going to give a
fairly tight lower bound on the worst case cost of evaluating a game tree
for the presence of a strategy.

The general case of determining the provability of a sequent can also be
organized into the question of whether a strategy exists in a formula,
although we have to construct a new game tree to represent this problem.
Given an opponent sequent $O \vd_o O'$, we show how to create a game $G$
(which is in fact $\ol{O} \otr O'$) so there is a proof of $O \vd_o O'$ if
and only if there is a strategy in $G$.

The translation between a proof of $O \vd_o O'$ and a strategy in $G$ is as
follows. Consider a proof for:

\begin{itemize}
\item $O \vd_o \bigsqcap_{i \in I} b_i \cn P_i$. It must start as
$(b_i \cn u_i)_{i \in I}$ where $u_i$ is a proof of $O \vd P_i$. If $I =
\{1,\ldots,m\}$ this may be represented as a strategy in the following tree:
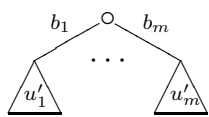
\[\xymatrix@R=2ex@C=4ex@M=0ex{
& \circ \ar@{-}[dl]_{b_1} \ar@{-}[dr]^{b_m} \\ 
\prooftri{\scriptstyle{u'_1}} & \cdots & \prooftri{\scriptstyle{u'_m}}}
\]
\medskip

where $u'_i$ is the tree whose strategies correspond to proofs of $O \vd P_i$.

\item $\bigsqcap_{i \in I} b_i \cn P_i \vd \bigsqcup_{j \in J} a_j \cn O_j$.
In this case we must make a move in either the left or right game and then 
find a strategy for the remaining problem. If $I = \{1,\ldots,m\}$ and $J =
\{1,\dots,n\}$ this may be represented as the following tree of choices:
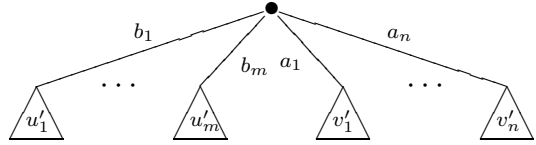
\[\xymatrix@M=0ex{
&&& \bullet \ar@{-}[dlll]_-{b_1} \ar@{-}[dl]^-{b_m} \ar@{-}[dr]_-{a_1}
\ar@{-}[drrr]^-{a_n} & \\ 
\prooftri{\scriptstyle{u'_1}} & \cdots & \prooftri{\scriptstyle{u'_m}} &&
\prooftri{\scriptstyle{v'_1}} & \cdots & \prooftri{\scriptstyle{v'_n}}}
\]
\medskip

where $u'_i$ is the tree whose strategies correspond to proofs of 
$P_i \vd \bigsqcup_{j \in J} a_j \cn O_j$ and
$v'_j$ is the tree whose strategies correspond to proofs of 
$\bigsqcap_{i \in I} b_i \cn P_i \vd O_j$.

\item $\bigsqcup_{j \in J} a_j \cn O_j \vd_p P$. It must start as
$\{a_j \cn v_j\}_{j \in J}$ where $v_j$ is a strategy for $O_j \vd P$. If
$J = \{1,\ldots,n\}$ this may be represented as the following tree:
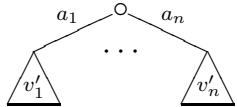
\[\xymatrix@R=2ex@C=5ex@M=0ex{
& \circ \ar@{-}[dl]_{a_1} \ar@{-}[dr]^{a_n} & \\
\prooftri{\scriptstyle{v'_1}} & \cdots & \prooftri{\scriptstyle{v'_n}}}
\]
\smallskip

where $v'_j$ is the tree whose strategies correspond to proofs of $O_j \vd P$.
\end{itemize}

This means the proof of any opponent sequent $O \vd_o O'$ can be transformed
into a strategy in $\ol{O} \otr O'$ using this translation process. From a
categorical (and logical) perspective this is just using the fact that:
\[\infer={(\,) \vd_o \ol{O} \otr O'}{O \vd_o O'}
\]
It would therefore seem that the complexity of provability may depend on
the size of this formula. This makes it interesting to know how the size of
this formula grows compared to the size of the original expressions.  

\section{On the size of multiplicatives} \label{sec-size}

To answer the question of how formulas grow when they are combined using
the multiplicative (and exponential) operators it is useful to develop some
techniques for calculating with game trees.

\subsection{Using multiplicities}

To calculate with trees a useful technique is to use multiplicities. 
Intuitively if a subgame occurs multiple times, rather than write the 
game multiple times, one can simply multiply it by the number of times 
it occurs. For example,
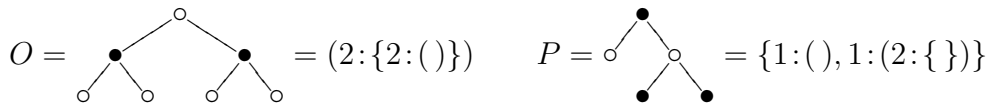
\[O =
\vcenter{\xymatrix@R=2ex@C=1.2ex@M=0ex{
&&& \circ \ar@{-}[dll] \ar@{-}[drr] \\
& \bullet \ar@{-}[dl] \ar@{-}[dr] &&&& \bullet \ar@{-}[dl] \ar@{-}[dr] \\
\circ && \circ && \circ && \circ}}
= (2 \cn \{2 \cn (\,)\}) \qquad P =
\vcenter{\xymatrix@R=2ex@C=1.2ex@M=0ex{
& \bullet \ar@{-}[dl] \ar@{-}[dr] \\
\circ && \circ \ar@{-}[dl] \ar@{-}[dr] \\
& \bullet && \bullet}}
= \{1 \cn (\,),1 \cn (2 \cn \{\,\})\}
\]
We may then calculate
\begin{align*}
O \oxr P &= \{1 \cn O \ox (\,),\, 1 \cn O \ox (2 \cn \{\,\})\} \\
&= \{1 \cn O,\, 1 \cn (2 \cn \{2 \cn (\,)\}) \ox (2 \cn \{\,\})\} \\
&= \{1 \cn O,\,1 \cn (2 \cn \{2 \cn (\,)\} \oxl (2 \cn \{\,\}),
   (2 \cn \{2 \cn (\,)\}) \oxr 2 \cn \{\,\})\} \\
&= \{1 \cn O,\,1 \cn (2 \cn \{2 \cn (\,) \ox (2 \cn \{\,\})\},2 \cn \{\,\})\} \\
&= \{1 \cn O,\, 1 \cn (2 \cn \{2 \cn (2 \cn \{\,\})\},2 \cn \{\,\})\} \\
&= \{1 \cn (2 \cn \{2 \cn (\,)\}), 1 \cn (2 \cn \{2 \cn (2 \cn \{\,\})\},
   2 \cn\{\,\})\}
\end{align*}
which, graphically, is the following tree:
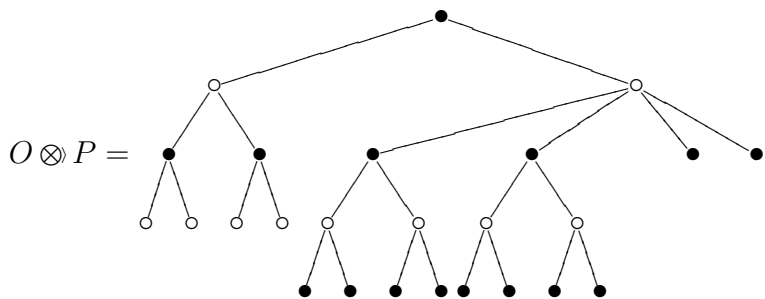
\[O \oxr P =
\vcenter{\xymatrix@R=4ex@C=0.5ex@M=0ex{
&&&&&&&&&&&&& \bullet \ar@{-}[dllllllllll] \ar@{-}[drrrrrrrrrr] \\
&&& \circ \ar@{-}[dll] \ar@{-}[drr] &&&&&&&&&&&&&&&&&&&& \circ
\ar@{-}[dlllllllllllll] \ar@{-}[dllllll] \ar@{-}[drrrrrr]
\ar@{-}[drrrrrrrrrrrrr] \\
& \bullet \ar@{-}[dl] \ar@{-}[dr] &&&& \bullet \ar@{-}[dl] \ar@{-}[dr]
&&&&& \bullet \ar@{-}[dll] \ar@{-}[drr]
&&&&&&& \bullet \ar@{-}[dll] \ar@{-}[drr]
&&&&&&&&&&&& \bullet &&&&&&& \bullet \\
\circ && \circ && \circ && \circ && \circ \ar@{-}[dl] \ar@{-}[dr]
&&&& \circ \ar@{-}[dl] \ar@{-}[dr] &&& \circ \ar@{-}[dl] \ar@{-}[dr]
&&&& \circ \ar@{-}[dl] \ar@{-}[dr] \\
&&&&&&& \bullet && \bullet && \bullet && \bullet & \bullet && \bullet &&
\bullet && \bullet}}
\]

While this does facilitate the calculation of multiplicatives, it is still
too detailed to allow an easy estimation of the size of a tensor or par of
two game trees. However, before abandoning this let us make an observation
on the exponentials.  Consider the following game:
\[!(n:\{m:(\,) \} )
\]
It is can be seen that this is the formula:
\[(n: \{m: (n-1: \{m: (n-2: \cdots (1:\{m:(\,)\}) \cdots )\})\})
\]
The leaf size (which is one more than the uniform size) of this formula is 
\[\leaves[!(n:\{m:(\,) \} )] = n! \cdot m^n
\]

This provides a lower bound on the possible growth in the size of a game tree
as this is a lower bound on the edge size of a particular tree.  

It is reasonable to wonder what an upper bound on the edge size might be. An
alternative way of viewing the exponential bang $(!)$ is as the game tree of
opponent initiated backtrackings. This means that on the opponents move one
can switch to any opponent move which had been a possibility on the play but
had not been chosen. This means a path is a sequence of choices of opponent
moves interleaved with a possible player move which may be played after the
opponent move. For $\esize_p[G] > 0$ there are a maximum of 
\[\esize_o[G]! \cdot \esize_p[G]^{\esize_o[G]}
\]
leaves, and hence paths. As each such path can contain at most
$2 \cdot \esize_o[G]$ edges the upper bound on the edge size is
\[2 \cdot \esize_o[G] \cdot \esize_o[G]! \cdot \esize_p[G]^{\esize_o[G]}
\]
When $\esize_p[G] = 0$ there are a maximum of $\esize_o[G]$ leaves which
is exactly the same as the number of edges in $G$. For games of depth two
or greater (i.e., $\esize_p[G] > 0$) however this shows that the
exponentials may produce very large game trees!

\subsection{Using profiles}

If $G$ is an arbitrary game tree its \textbf{profile} is a list of the number
of nodes which occur at each level. Notice the list always begins with a $1$.
The two trees above and their tensor have the following profiles:
\[
\profile[O] = [1,2,4], \quad \profile[P] = [1,2,2],\quad \text{and} \quad
\profile[O \oxr P] = [1,2,6,8,8]
\]
The notation $\profile[G]_i$ will be used to indicate the number of nodes
of $G$ at the $i$-th depth (starting at zero). We have the following useful
observation:

\begin{proposition} For $O$ and $O'$ opponent games and $P$ a player game

\begin{enumerate}[{\upshape (i)}]
\item $\profile[O \ox O']_0 = \profile[O \otr P]_0 = 1$

\item $\profile[O \ox O']_n = \sum\limits^n_{i=0}
\choose {\lfloor n/2 \rfloor} {\lfloor i/2 \rfloor} \cdot
\gamma^i_{n-i} \cdot \profile[O]_i \cdot
\profile[O']_{n-i}$ 

\item $\profile[O \oxr P]_{n+1} = \sum\limits^n_{i=0} 
\choose{\lfloor n/2 \rfloor}{\lfloor i/2 \rfloor} \cdot
\gamma^i_{n-i} \cdot \profile[O]_{i} \cdot \profile[P]_{n-i+1}$
\end{enumerate}
where 
\[\gamma_n^m = 
\begin{cases}
1 & \text{if $n$ or $m$ is zero or even} \\
0 & \text{otherwise}
\end{cases}
\]
\end{proposition}

The rationale for this formula is as follows. To get a node at depth $n$ in
$O \ox O'$ we must interleave the path to a node of $O$ at depth $i$
with the path to a node of $O'$ at depth $n-i$. As we can interleave only
at the even nodes there are $\choose{\lfloor n/2 \rfloor}{\lfloor i/2 \rfloor}$
possible such interleavings all of which provide distinct nodes. Finally, we
can only produce interleaved paths of the required length
if at least one of $i$ or $n-i$ is even as we have to end in one of the
paths having left the other at an even index. This results in the factor
$\gamma^n_i$ which is applied to each term.

\begin{proof}
It is easy to see that the formulas work for the base cases $(\,)$ and
$\{\, \}$. Let $O = (a_i:P_i \mid 0 \leq i < n)$ and $O' = (b_j:P'_j \mid 0
\leq j < m)$ then 
\[O \ox O' = (a_i:P_i \oxl O', b_j:O \oxr P'_j \mid 0 \leq i < n, 0 \leq j < m)
\]
For depths of zero and one its profile is respectively
\[\profile[O \ox O']_0  = 1 \quad \text{and} \quad \profile[O \ox O']_1 = n+m
\]
Checking that the formula is correct for these depths is straightforward.
Now consider the inductive step for the last case:
\begin{equation*}
\begin{split}
\profile[O \ox O']_{r+2}
&= \sum_{i=1}^{n} \profile[P_i \oxl O']_r + \sum_{j=1}^{m}
    \profile[O \oxr P'_j]_r \\
&=\sum_{k=1}^{n} \sum^r_{i=0}
    \choose{\lfloor \frac{r}{2} \rfloor}{\lfloor \frac{i}{2} \rfloor} \cdot
    \gamma_i^{r-i} \cdot \profile[O']_i \cdot \profile[P_k]_{r-i+1} \\
&\quad  + \sum_{j=1}^{m} \sum^r_{i=0} 
    \choose{\lfloor \frac{r}{2} \rfloor}{\lfloor \frac{i}{2} \rfloor} \cdot
    \gamma_i^{r-i} \cdot \profile[O]_i \cdot  \profile[P_j]_{r-i+1} \\
&= \sum^r_{i=0}
    \choose{\lfloor \frac{r}{2} \rfloor}{\lfloor \frac{i}{2} \rfloor} \cdot
    \gamma_i^{r-i} \cdot \profile[O']_i \cdot
    \left(\sum_{k=1}^{n} \profile[P_k]_{r-i+1}\right) \\
&\quad + \sum^r_{i=0} 
    \choose{\lfloor \frac{r}{2} \rfloor}{\lfloor \frac{i}{2} \rfloor} \cdot
    \gamma_i^{r-i} \cdot \profile[O]_i \cdot
    \left(\sum_{j=1}^{m} \profile[P'_j]_{r-i+1}\right) \\
&= \sum^r_{i=0} 
    \choose{\lfloor \frac{r}{2} \rfloor}{\lfloor \frac{i}{2} \rfloor} \cdot
    \gamma_i^{r-i} \cdot \profile[O']_i \cdot \profile[O]_{r-i+2} \\
&\quad + \sum^r_{i=0} 
    \choose{\lfloor\frac{r}{2}\rfloor}
           {\lfloor\frac{r}{2}\rfloor-\lfloor\frac{i}{2}\rfloor} \cdot
    \gamma_i^{r-i} \cdot \profile[O]_i \cdot \profile[O']_{r-i+2} \\
\intertext{as $\lfloor \frac{r}{2} \rfloor - \lfloor \frac{i}{2} \rfloor = 
\lfloor \frac{r-i}{2} \rfloor$ we may reindex the second sum:}
&= \sum^r_{i=0} 
    \choose{\lfloor \frac{r}{2} \rfloor}{\lfloor \frac{i}{2} \rfloor} \cdot
    \gamma_i^{r-i+2} \cdot \profile[O']_i \cdot \profile[O]_{r-i+2} \\
&\quad + \sum^{r+2}_{i'=2} 
    \choose{\lfloor\frac{r}{2}\rfloor}{\lfloor\frac{i'}{2}\rfloor-1} \cdot
    \gamma_{i'}^{r-i'+2} \cdot \profile[O']_{i'} \cdot \profile[O]_{r-i'+2} \\
&= \sum^{r}_{i=2}
    \left(\choose{\lfloor \frac{r}{2} \rfloor}{\lfloor \frac{i}{2} \rfloor-1} +
    \choose{\lfloor \frac{r}{2} \rfloor}{\lfloor \frac{i}{2} \rfloor}\right)
\cdot \gamma_i^{r-i+2} \cdot \profile[O]_{i} \cdot \profile[O']_{r-i+2} \\
&\quad  + \choose{\lfloor\frac{r}{2}\rfloor}{0} \cdot
    (\gamma_0^{r+2} \cdot \profile[O']_0 \cdot \profile[O]_{r+2} 
    + \gamma_1^{r+1} \cdot \profile[O']_1 \cdot \profile[O]_{r+1})  \\
&  + \choose{\lfloor\frac{r}{2}\rfloor}{\lfloor\frac{r}{2}\rfloor} \cdot
      (\gamma^1_{r+1} \cdot \profile[O']_{r+1} \cdot \profile[O]_1 
      \gamma^1_{r+1} \cdot \profile[O']_{r+2} \cdot \profile[O]_0) \\
&= \sum^{r+2}_{i=0}
    \choose{\lfloor \frac{r}{2} \rfloor+1}{\lfloor \frac{i}{2} \rfloor} \cdot
    \gamma_i^{r-i+2} \cdot \profile[O]_{i} \cdot \profile[O']_{r-i+2} 
\end{split}
\end{equation*}
\end{proof}

We now wish to determine how the cost of evaluating tensors of trees can
grow. A first guess might be that the tensor always causes an exponential
increase in size. However, this is not the case in general: the growth 
depends on the shape of the tree. Below we consider two examples where an
exact computation of the uniform size (hence our interest in this measure)
is possible:

\begin{enumerate}
\item First we consider a case where the growth is polynomial:
\[A_0 = \{\,\},\ A_{n+1} = \{2 \cn B_n\} \quad \text{and} \quad
B_0 = (\,),\ B_{n+1} = (1 \cn A_n) 
\]
so that, for example,
\[A_3 = \{2 \cn (1 \cn \{2 \cn (\,)\})\}\ =
\vcenter{\xymatrix@R=2ex@C=0.7ex@M=0ex{
&&&&&&& \bullet \ar@{-}[dllll] \ar@{-}[drrrr] \\
&&& \circ \ar@{-}[d]  &&&&&&&& \circ \ar@{-}[d] \\
&&& \bullet \ar@{-}[dl] \ar@{-}[dr] &&&&&&&& \bullet \ar@{-}[dl] \ar@{-}[dr] && \\
 && \circ && \circ &&&&&& \circ && \circ &&}}
\]

Notice that $\usize[A_{2n}] = \usize[A_{2n-1}]$ and so $\usize[A_{2n}] =
2^n-1$. Now consider $\usize[A_{2n} \ot A_{2m}]$. The tree $A_{2n} \ot
A_{2m}$ has a uniform maximum depth of $2(n+m)$ whose leaf count is
concentrated in the term 
\begin{align*}
\choose{n+m}{n} \cdot \gamma^n_m \cdot \profile[A_{2n}]_{2n} \cdot
\profile[A_{2n}]_{2m}
&= \choose{n+m}{n} \cdot 2^n \cdot 2^m \\
&\leq 2^{n+m} \cdot 2^n \cdot 2^m \\
&=  (2^n)^2 \cdot (2^m)^2
\end{align*}
Thus, the uniform size is bounded by a (degree 4) polynomial in the uniform
sizes of the two games in this case.

\item To demonstrate that an exponential growth is possible we consider a
very simple example:
\[L_0 = \{\,\},\ L_{n+1} = \{1 \cn L'_n\} \quad \text{and} \quad
L'_0 = (\,),\ L'_{n+1} = (1 \cn L_n) 
\]

This gives $|L_{2n}| = |L_{2n}|_p + |L_{2n}|_o = (n+1) + n$ although
$\usize[L_{2n}] = 0$. However
\begin{align*}
\usize[L_{2n} \ot L_{2m}] &= \profile[L_{2n} \ot L_{2m}]_{2n+2m} \\
&= \choose{n+m}{n} \cdot \profile[L_{2n}]_{2n} \cdot \profile[L_{2m}]_{2m} \\
&= \choose{n+m}{n} \\
&\geq \frac{n^m + m^n}{2} 
\end{align*}

This gives a lower bound on the uniform size of the tensored trees which
has an exponential increase on the number of nodes in the trees.
\end{enumerate}

This gives an indication that the complexity of a naive evaluation can be
quite bad.

\section{Dynamic programming and graph games} \label{sec-dyn}

We have already established that to determine whether a sequent $(\,) \vd_o
U$, where $U$ is an arbitrary additive formula, has a proof may be
determined in time bounded by (and in the worst case) $\esize[U]$.
Fortunately, for the provability of $O \vd_o
O'$, where $O$ and $O'$ are opponent games, there is something further we
can do as the structure of $O \vd_o O'$ can be exploited.

\subsection{Dynamic programming}

Notice in the calculation of $\ol{O} \otr O'$ that subexpressions may
occur multiple times. This indicates that it may be valuable to attempt a
dynamic programming approach where one calculates the values of the
sub-problems at most once. This is, in fact, the case and in this section
we show:

\begin{proposition} \label{prop-dyn}
Using dynamic programming one can evaluate the provability of a sequent
$O \vd_o O'$ using
\[\usize[O] |O'|_p + \usize[O'] |O|_p + |O|_p |O'|_p
\]
binary operations.
\end{proposition}

Notice that this is a substantial improvement on a simple evaluation of the 
game tree resulting from the internal-hom object. The calculation below
essentially counts the number of binary operations required but is generous
enough to allow for single argument evaluations. An initial observation is:

\begin{lemma} \label{lemma1}
For any player game 
$P = \{a_i \cn (b_j \cn P_{ij} \mid j \in \{1,\ldots,n_i\})
\mid i \in \{1,\ldots,n\}\}$:

\begin{enumerate}[{\upshape (i)}]
\item $\sum\limits_{i = 1}^n \sum\limits_{j = 1}^{n_i} |P_{ij}|_p = |P|_p - 1$
\item $\usize[P] = n-1 + \sum\limits_{i=1}^n (n_i-1) +
\sum\limits_{i = 1}^n \sum\limits_{j = 1}^{n_i} \usize\left[P_{ij}\right]$
\end{enumerate}
\end{lemma}

\begin{proof}[of Proposition~\ref{prop-dyn}]
To determine the evaluation of $P \otr O$, where $O =
(a_i: P_i \mid i \in \{1,\ldots,n\})$, we distribute $P$ inside the opponent
game to obtain $(a_i:P \ot P_i \mid i \in \{1,\ldots,n\})$. This can be
solved in $n-1 + \sum_{i=1}^n \cost[P \ot P_i]$ logical operations, where
$\cost[G]$ denotes the cost to evaluate $\|G\|$. In what follows we show
\[\cost[P \ot P_i] = \usize[P] |P_i|_p + \usize[P_i] |P|_p + |P|_p |P_i|_p
\]
so that
\begin{align*}
\cost[P \otr O] &= \cost[P \otr (a_i: P_i \mid i \in \{1,\ldots,n\})] \\
&= n-1 + \usize[P] \cdot \sum_{i=1}^n |P_i|_i + |P|_p \cdot \sum_{i=1}^n
\usize[P_i]_i + |P|_p \cdot \sum_{i=1}^n |P_i|_i \\
&= n-1+ \usize[P] |O|_p + |P|_p (\usize[O]-(n-1)) + |P|_p |O|_p \\
&\leq \usize[P] |O| + \usize[O] |P|_p + |P|_p |O|_p \\
\end{align*}

For player games $P$ and $P'$ we now calculate $\cost[P \ot P']$. First fix
some notation:
\[\begin{array}{ccc}
P = \{a_i \cn O_i \mid i \in \{1,\ldots,m\}\} & \quad &
P' = \{b_j \cn O'_j \mid j \in \{1,\ldots,n\}\}
\medskip\\
O_i = (a_{ij} \cn P_{ij} \mid j \in \{1,\ldots,m_i\})
& \quad &
O'_j = (b_{ji} \cn P_{ji} \mid i \in \{1,\ldots,n_j\})
\end{array}
\]

The calculation will proceed as follows:
\begin{align*}
\| P \ot P' \| &= \Big(\bigjoin_{i=1}^m \| O_i \otl P' \| \Big) 
\join \Big(\bigjoin_{j=1}^n \| P \otr O_j' \| \Big) \\
&= \Big(\bigjoin_{i=1}^m \Big(\bigmeet_{k=1}^{m_i}\| P_{ik} \ot P'\|\Big)\Big)
\join 
\Big(\bigjoin_{j=1}^n \Big(\bigmeet_{l=1}^{n'_j} \| P \ot P_{jl}'\|\Big)\Big)
\end{align*}

This shows how this value can be calculated using dynamic programming
techniques since the computation of $\|P_{ik} \ot P'_{jl}\|$ will reoccur.

We will now calculate inductively the number of logical operations involved
in calculating this value. Note that the base cases are when the index sets
are empty; these are handled by the general argument. The following ``tree of
required operations'' may be helpful to visualize the calculation.

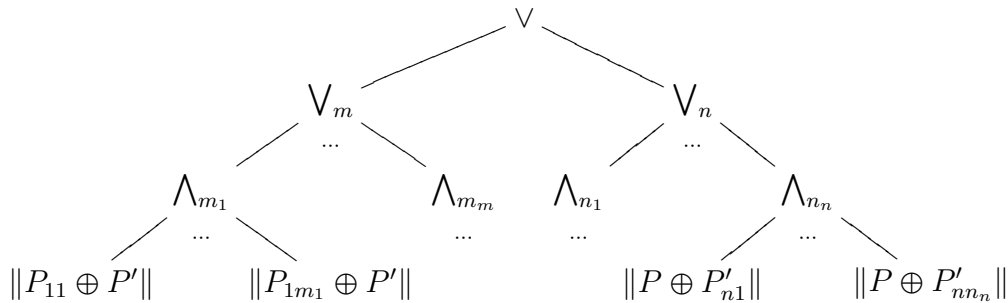
\[\xymatrix@R=3ex@C=0.1ex{
&&&& \join \ar@{-}[dll] \ar@{-}[drr] \\
&& \bigjoin_m \ar@{-}[dl] \ar@{-}[dr] \ar@{}[d]|\cdots 
&&&& \bigjoin_n \ar@{-}[dl] \ar@{-}[dr] \ar@{}[d]|\cdots  \\
& \bigmeet_{m_1} \ar@{-}[dl] \ar@{}[d]|\cdots \ar@{-}[dr] 
&& \bigmeet_{m_m} \ar@{}[d]|\cdots
&& \bigmeet_{n_1} \ar@{}[d]|\cdots
&& \bigmeet_{n_{n}} \ar@{-}[dl] \ar@{}[d]|\cdots \ar@{-}[dr] \\
\| P_{11} \ot P' \| && \| P_{1m_1} \ot P' \| &&&& \| P \ot P'_{n1} \|
&& \| P \ot P'_{nn_{n}} \|}
\]

It takes $m-1$ logical operations to form the join or meet of $m$ terms. 
The calculations for the left side of the tree
\[(m-1) + \sum_{i=1}^m (m_i-1) + \sum_{i=1}^m \sum_{k=1}^{m_i}
\cost[P_{i_k} \ot P']
\]
plus the calculations for the right side of the tree
\[(n-1) + \sum_{j=1}^n (n_j-1) + \sum_{j=1}^n \sum_{l=1}^{n_j}
\cost[P \ot P_{j_l}']
\]
plus the one join at the root, therefore, gives us the total number of
logical operations. Notice that if either or both sides of this calculation
reduce to a single argument calculation we still have the cost of the
single binary operation into which the retrieval costs can be bundled.

In the calculations of the $\|P_{ik} \ot P'\|$'s and the $\|P \ot P_{jl}'\|$'s
it can be seen that we will end up calculating the $\|P_{ik} \ot P_{jl}'\|$'s
twice. Thus, from the overall sum we can subtract the number of calculations
required to calculate this value once, which is
\[\sum_{i=1}^m\sum_{k=1}^{m_i} \sum_{j=1}^n\sum_{l=1}^{n_j}
\cost[P_{ik} \ot P'_{jl}]
\]
What we then have is:
\begin{align*}
\cost[P \ot P'] &= 1+(m-1)+ \sum_i(m_i-1) + \sum_{ik} \cost[P_{ik} \ot P'] \\
&\quad + (n-1) + \sum_j (n_j-1) + \sum_{jl} \cost[P \ot P'_{jl}] \\
&\quad -\sum_{ik} \sum_{jl} \cost[P_{ik} \ot P'_{jl}]
\end{align*}
Applying the theorem inductively to the three embedded terms yields
\begin{align*}
&\leq 1 + (m-1) + \sum_i (m_i-1) + \sum_{ik} (\usize[P_{ik}]  |P'|_p +
\usize[P']  |P_{ik}|_p  + |P_{ik}|_p  |P'|_p) \\
&\quad + (n-1) + \sum_j (n_j-1) + \sum_{jl} (\usize[P]  |P_{jl}'|_p +
\usize[P_{jl}']  |P|_p  + |P|_p  |P'_{jl}|_p) \\
&\quad - \sum_{ik} \sum_{jl} (\usize[P_{ik}]  |P'_{jl}|_p+ \usize[P'_{jl}] 
 |P_{ik}|_p  + |P_{ik}|_p  |P'_{ij}|_p)
\end{align*}
Applying Lemma~\ref{lemma1} (i) and rearranging terms: 
\begin{align*}
&= 1 + (m-1) + \sum_i(m_i-1) \\
&\quad + \left(\sum_{ik} \usize[P_{ik}] |P'|_p - \usize[P_{ik}] (|P'|_p-1)
       \right) + \usize[P]  (|P'|_p-1) ) \\
&\quad + (n-1) + \sum_i(n_i-1) \\
&\quad + \left(\sum_{jl} \usize[P'_{jl}] |P|_p - \usize[P'_{jl}] (|P|_p-1)
       \right) + \usize[P']  (|P|_p-1) \\
&\quad + |P|_p |P'|_p
\end{align*}
With some more rearranging:
\begin{align*}
&= (m-1) + \sum_i(m_i-1) +\sum_{ik} \usize[P_{ik}] + \usize[P]  (|P'|_p-1) \\
&\quad +(n-1)+\sum_i(n_i-1)+\sum_{jl} \usize[P'_{jl}] + \usize[P'] (|P|_p-1) \\
&\quad + |P|_p  |P'|_p
\end{align*}
And finally, by Lemma~\ref{lemma1} (ii):
\[
\leq \usize[P] |P'|_p + \usize[P'] |P|_p  + |P|_p |P'|_p
\]
\end{proof} 

Usually to calculate the complexity of a dynamic programming algorithm one
counts the number of subproblems and multiplies that number by the cost of
computing that problem. As the cost of computation at each node varies
according to the structure of the games we have done a more detailed
analysis which removes the cost of the common computation. We now return to
apply these counting techniques with confidence!

\subsection{Graph games}

This dynamic programming approach encourages us to view the games, not as
trees, but as acyclic directed graphs. This is, of course, not a new idea
and has been explored by Hyland and Schalk \cite{hyland-schalk}. The point
is that one of these ``graph games'' will have a strategy if and only if
the corresponding ``tree game'' has such. To evaluate whether a graph game
has a strategy requires time proportional to the number of edges in the
graph.

One of the attractions of graph games is that one can count the number of
edges and nodes of any game which is expressed using the multiplicatives
(exponentials seem harder to analyze as they use the history of how one
reaches a node). This means that one can calculate the cost of evaluating
(on this initial model) an arbitrary polarized multiplicative and additive
formula using dynamic programming techniques surprisingly easily. 

We should perhaps remind the reader at this stage that the multiplicative
structure of these models is emphatically not free. Provability in the free
such structure with multiplicative units is NP-complete as this is the case
in the unpolarized setting (due to Kanovich's
results~\cite{kanovich92:complexity}). Recall also that with exponentials
and additives in the unpolarized setting these problems become undecidable.
It is likely, although we have not done the calculation, that this is also
true for the free polarized multiplicative-additive logics with exponentials.

In contrast to the above results the problem in the initial polarized game
logic is, for all these fragments, quite decidable as a given formula can
be ``unwrapped'' into its purely additive structure. The problem is that
unwrapping into purely additive structure can result in a huge increase in
the size of the formula. Given the analysis above it is reasonable to
suppose one can do much better then naively expanding into a game tree and
evaluating.

We shall now expect our formula to belong to the polarized additive,
multiplicative fragment. In other words to be in the following form:

\begin{center}
\begin{tabular}{c|c}
Opponent & Player \\ \hline & \\
$( b_i: P_i \mid i \in I) \qquad O \ox O'$ &
$\{a_j: O_j \mid j \in J \} \qquad P \ot P'$ \\ \quad & \quad \\
$O \otl P \qquad P \otr O$ & $P \oxl O \qquad O \oxr P$ \\
\end{tabular}
\end{center}

It is easy to see that the cost of evaluation in a graph game is dominated
by the number of edges in the game. Thus, for example, in the graph game
produced for the dynamic calculation of $P \ot P'$ there are $|P|_p |P'|_p$
player states. From each of these states one can either make a move in the
first or second coordinate resulting in
\[\esize_p[P \ot P] = \esize_p[P]  |P'|_p + \esize_p[P'] |P|_p
\]
We use the fact that the number of edges leaving a node (up to a factor)
always dominates the cost of the calculation needed at that node. Similarly
there are $|P|_o |P'|_p + |P|_p |P'|_o$ opponent states. While the number
of opponent edges is given by:
\[\esize_o[P \ot P] = \esize_o[P] |P'|_p + \esize_o[P'] |P_1|_p
\]
We therefore have:

\begin{lemma}
The cost of evaluating multiplicative expressions using dynamic programming 
(up to a constant factor) is bounded by:

\begin{itemize}
\item $\cost[(b_i: P_j \mid j \in \{1,..,m\})] = m+\sum^m_{j=1} \cost[P_j]$

\item $\cost[\{a_i: O_i \mid i \in \{1,..,n\})] = n+\sum^n_{i=1} \cost[O_i]$

\item $\cost[O \ox O] = \esize[O] |O'|_o + \esize[O'] |O|_o$

\item $\cost[O \oxr P] = \esize[O] |P|_o + \esize[P] |O|_o$

\item $\cost[P \ot P'] = \esize[P] |P'|_p + \esize[P'] |P|_p$

\item $\cost[P \otr O] = \esize[P] |O|_p + \esize[O] |P|_p$
\end{itemize}
\end{lemma}

Note that in this calculation there is no need to assume that we are
starting with trees. We may start with acyclic directed graphs and
recursively apply the calculation. Indeed, to evaluate an arbitrary
expression this is what we must do: a multiplicative expression allows us
to build an acyclic graph with a certain number of player and opponent
nodes and arrows. As we apply constructors we build new graphs, all we
need to keep track of is the numbers of opponent nodes, player nodes and
edges. Below is the explicit calculation for doing so. Let 
\[\size[G] = (|G|_o,|G|_p,\esize_p[G])
\]
and $\size[G_1] = (n_1,N_1,E_1)$ and $\size[G_2] = (n_2,N_2,E_2)$.
Then:

\begin{itemize}
\item $\size[(b_1:G_1,\ldots,b_n:G_n)] = \sum^n_{i=1} \size[G_i] + (0,1,n)$

\item $\size[\{a_1:G_1,\ldots,a_m:G_m\}]= \sum^m_{j=1} \size[G_j] + (1,0,m)$

\item $\size[G_1 \ot/\otl/\otr G_2] = (n_1 \cdot n_2,\
n_1 \cdot N_2 + N_1 \cdot n_2,\ E_1 \cdot n_2 + n_1 \cdot E_2)$

\item $\size[G_1 \ox/\oxl/\oxr G_2] = (N_1 \cdot n_2 + n_1 \cdot N_2,\
N_1 \cdot N_2,\ E_1 \cdot N_2 + N_1 \cdot E_2)$
\end{itemize}

The cost of determining whether a formula is provable is bounded by the last
entry in this cost calculation. It is worth noting that, even though the
dynamic approach to this evaluation has helped, we are still left with an
exponential bound for the evaluation of the provability of a multiplicative
expression as the complexity for multiplicatives is essentially multiplicative!

We may also very easily calculate the profile of the graph games which results
from multiplicative expressions. The only change needed over our previous
calculation is to remove the binomial coefficients. For profiles the only
distinction between player and opponent is the starting point of the profile.
It is necessary, therefore, only to state the calculation for tensor and sums:


\begin{itemize}
\item $\profile[\{a_1: G_1,...,a_n: G_n \}]_0 = 1$
\item $\profile[\{a_1:G_1,...,a_n:G_n\}]_{m+1} = \sum^n_{i=1} \profile[G_i]_m$
\item $\profile[G_1 \ox G_2]_n = \sum^n_{i=0} \gamma^i_{n-i} \cdot
       \profile[G_1]_i \cdot  \profile[G_2]_{n-i}$
\item $\profile[G_1 \oxr P]_{0} = 1$
\item $\profile[G_1 \oxr P]_{n+1} = \sum^n_{i=0} \gamma^i_{n-i} \cdot
\profile[G_1]_{i} \cdot  \profile[P]_{n-i+1}$
\end{itemize}

This allows an easy calculation of the number of states used by the dynamic
algorithm. However, notice that the profile does lose the edge information.
One can actually bound the number of edges  by:
\[\esize(G) \leq \sum^{n-1}_{i=0} \profile[G]_i \cdot \profile[G]_{i+1}
\]
as each edge must end in a level one greater than where it starts. Of course
this bound is going to be, in general, very poor; its only advantage is that
it can be derived from the profile of the acyclic graph which is easily
calculated.

\section{The complexity of provability}  \label{sec-prov}

The technique of viewing a game tree as an acyclic graph has certainly led
to improvements in the determination of provability. We hope also it has
led to the reader having a better feel for the combinatoric structure of
these games. 

However, the story of the complexity of provability is not over! We have
yet to really use the categorical structure the initial model of polarized
games which is what we now exploit. We open with the following observation:

\begin{proposition}
In the initial model of polarized games:

\begin{enumerate}[{\upshape (i)}]
\item A game $G$ has a strategy if and only if there is (a possibly mixed)
map $(\,) \to G$.
\item A game $G$ has a counter-strategy if and only if there is a (a
possibly mixed) map $G \to \{\, \}$.
\item Every game either has a strategy or a counter-strategy but not both.
\item $O \ox O'$ has a strategy if and only if $O$ and $O'$ have
strategies. Similarly, $O \oxr P$ has a strategy if and only if $O$ and $P$
have strategies.
\item $P \ot P'$ has a counter-strategy if and only if $P$ and $P'$ have
counter-strategies. Similarly, $P \otr O$ has a counter-strategy if and
only if $O$ and $P$ have counter-strategies.
\item $\ol{G}$ has a strategy (resp. counter-strategy) if and only if $G$
has a counter-strategy (resp. strategy).
\item $!O$ has a strategy if and only if $O$ has a strategy.
\item $?P$ has a counter-strategy if and only if $P$ has a counter strategy.
\end{enumerate}
\end{proposition}

\begin{proof} \quad
\begin{enumerate}[{\upshape (i)}]
\item This follows immediately from the discussion in Section~\ref{sec-strat}.

\item Dual to {\upshape (i)}.

\item If a game $G$ has a strategy and a counter-strategy then there would
be a mix map $(\,) \to \{\, \}$ which there is not! In Remark~\ref{rem-strat}
we noted that the calculation to determine the possession of a strategy is
exactly the de Morgan dual calculation to that for having a
counter-strategy. Thus, if one is true the other is necessarily false.

\item If $(\,) \to^s O \ox O'$ then, as $(\,)$ is a final object, we obtain 
\[(\,) \to^s O \ox O' \to^{1 \ox \mathbf{!}} O \ox (\,) = O
\]
and so $O$ and (similarly) $O'$ have strategies.  Conversely if $(\,)
\to^{s_1} O$ and $(\,) \to^{s_2} O'$ then 
\[(\,) = (\,) \ox (\,) \to^{s_1 \ox s_2} O \ox O'
\]
provides a strategy for $O_1 \ox O_2$. A similar argument works for the
mixed tensor: if $(\,) \to^{s_1} O$ and $(\,) \to^{s_2} P$ then
\[(\,) = (\,) \ox (\,) \to^{s_1 \ox s_2} O \ox P
\]
Conversely, if $(\,) \to^{s} O \oxr P$ then we may compose with $O \oxr P
\to^{\mathbf{!} \oxr 1_P} (\,) \oxr P = P$ to obtain a strategy in $P$. As
$P$ has a strategy $P \neq \{\,\}$ and so $P = \{a_i:O_i \mid i \in I\}$
where $I \neq \emptyset$. This means that in $O \oxr P = \{a_i:O \ox O_i \mid
i \in I\}$ the strategy $s$ must choose a branch and so there is a strategy
in $O \ox O_i$ for some $i \in I$. By the argument above this happens if and
only if there is a strategy in both $O$ and $O_i$. Thus, there is certainly
a strategy in $O$.

\item Dual to {\upshape (iv)}.

\item If $(\,) \to^s O$ then we obtain a strategy 
\[(\,)= !(\,) \to^{!s} !O
\]
Conversely, if $(\,) \to^{s'} !O$ is a strategy we can use dereliction to
obtain a strategy
\[(\,) \to^{s'} !O \to^{\epsilon} O
\]

\item Dual to {\upshape (v)}.
\end{enumerate}
\end{proof}

This means that we may completely redo the our calculation of the cost of
provability! To determine whether $O \vd_o O'$ is provable we must
determine whether there is a strategy in $\ol{O} \otr O'$. But this is
the case only if there is no counter-strategy for $\ol{O} \otr O'$.
There is a counter-strategy for $\ol{O} \otr O'$ only if there is a
counter-strategy for $\ol{O}$ and $O'$, that is there is no strategy for 
$O$ and there is a counter-strategy for $O'$. In other words, $O \vd_o O'$
is provable if and only if either there is no strategy for $O$ or a
strategy for $O'$.

Clearly we may now show that the cost of provability is \emph{linear} for
an arbitrary sequent in the polarized additive, multiplicative, and
exponential fragment by successively using the reductions provided by the
proposition. Therefore we have:

\begin{corollary}
In initial polarized game logic the complexity of provability is linear in
the size of the sequents.
\end{corollary}

\begin{proof}
We shall explicitly provide a procedure based on the proposition for
determining the provability of a formula. Let $\|P\|_p$ be the value of
provability of $P \vd_p \{\,\}$ and $\|O\|_o$ the value of provability of
$(\,) \vd_o O$. Then:
\begin{eqnarray*}
\|P \ot P'\|_p &=& \|P\|_p \join \|P'\|_p \\
\|O \oxr P\|_p = \|P \oxl O\|_p &=& \|O\|_o \meet \|P\|_p \\
\|\ol{O}\|_p &=& \neg \|O\|_o \\
\|\{a_i:O_i \mid i \in I\}\|_p &=& \bigjoin_{i \in I} \|O_i\|_o \\
\|O \ox O'\|_o &=& \|O\|_o \meet \|O'\|_o \\
\|O \oxl P\|_o = \|P \oxr O\|_o &=& \|O\|_o \join \|P\|_p \\
\|\ol{P}\|_o &=& \neg \|P\|_p \\
\|(b_j:P_j \mid j \in J\}\|_o &=& \bigmeet_{j \in J} \|P_j\|_p
\end{eqnarray*}
There are dual rules for determining whether there is a counter-strategy.
Clearly this procedure is linear in the size of the term.
\end{proof}

Rather surprisingly provability has therefore reduced to the cost of looking
at the formulas. As, in the worst case, this is the very least one must do
this is a hard bound.  



\end{document}